\documentclass{amsart}

\usepackage{a4wide,amsmath,amssymb}
\usepackage[toc,page]{appendix}

\usepackage[integrals]{wasysym}

\usepackage{url}

\usepackage{hyperref}

\usepackage{mathtools}

\usepackage{xcolor}

\newtheorem{thm}{Theorem}[section]
\newtheorem{lemma}[thm]{Lemma}

\newtheorem{cor}[thm]{Corollary}
\theoremstyle{remark}
\newtheorem{rem}[thm]{Remark}

\theoremstyle{definition}

\newtheoremstyle{Claim}{}{}{\itshape}{}{\itshape\bfseries}{:}{ }{#1}
\theoremstyle{Claim}

\newcommand{\T}{{\mathbb{T}}}

\newcommand{\R}{\mathbb{R}}

\newcommand{\eps}{\varepsilon}

\newcommand{\dd}{\hspace{0.7pt}{\rm d}}

\makeatletter
\@namedef{subjclassname@2020}{\textup{2020} Mathematics Subject Classification}
\makeatother

\theoremstyle{plain}

\begin{document}

\title[]{Convergence rates for the vanishing viscosity approximation of Hamilton-Jacobi equations:\\ the convex case}

\author{Marco Cirant}
\address{Dipartimento di Matematica ``Tullio Levi-Civita'', Universit\`a degli Studi di Padova, 
via Trieste 63, 35121 Padova (Italy)}
\curraddr{}
\email{cirant@math.unipd.it}
\thanks{}
\author{Alessandro Goffi}
\address{Dipartimento di Matematica e Informatica ``Ulisse Dini'', Universit\`a degli Studi di Firenze, 
viale G. Morgagni 67/a, 50134 Firenze (Italy)}
\curraddr{}
\email{alessandro.goffi@unifi.it}

\date{\today}


\begin{abstract}
We study the speed of convergence in $L^\infty$ norm of the vanishing viscosity process for Hamilton-Jacobi equations with uniformly or strictly convex Hamiltonian terms with superquadratic behavior. Our analysis boosts previous findings on the rate of convergence for this procedure in $L^p$ norms, showing rates in sup-norm of order $\mathcal{O}(\eps^\beta)$, $\beta\in(1/2,1)$, or $\mathcal{O}(\eps|\log\eps|)$ with respect to the vanishing viscosity parameter $\eps$, depending on the regularity of the initial datum of the problem and convexity properties of the Hamiltonian. Our proofs are based on integral methods and avoid the use of techniques based on stochastic control 
or the maximum principle.
\end{abstract}

\maketitle

\section{Introduction}

This paper studies quantitative estimates of the rate of convergence of the vanishing viscosity approximation of periodic solutions to the viscous Hamilton-Jacobi equation
\begin{equation}\label{HJintroVisc}
\begin{cases}
-\partial_t u_\eps-\eps\Delta u_\eps+H(Du_\eps)=f(x,t)&\text{ in }Q_T:=\T^n\times(0,T),\\
u_\eps(x,T)=u_T(x)&\text{ in }\T^n.
\end{cases}
\end{equation}
towards its first-order counterpart
\begin{equation}\label{HJintro}
\begin{cases}
-\partial_t u+H(Du)=f(x,t)&\text{ in }Q_T:=\T^n\times(0,T),\\
u(x,T)=u_T(x)&\text{ in }\T^n.
\end{cases}
\end{equation}
Here, $H:\R^n\to\R$ is the so-called Hamiltonian, while $f:Q_T\to\R$ and $u_T:\T^n\to\R$ are the periodic source and the terminal condition of the problem.\\

The rate of convergence of this regularization approach has been widely investigated in the context of viscosity solutions and it is sensitive of the regularity and the geometric properties of the Hamiltonian and the data of the problem. In the case of Lipschitz solutions and \textit{locally Lipschitz} $H$ the $\mathcal{O}(\sqrt{\eps})$ rate of convergence for $\|u_\eps-u\|_\infty$ was established in \cite{FlemingPrinceton,FlemingJMM,L82Book,CL84,EvansARMA,IshiiKoike} using probabilistic techniques, maximum principle approaches and integral duality methods respectively, see also \cite{Calder} for a proof via the regularization by sup/inf-convolution. The Lipschitz continuity of $H$ can be flattened to the mere continuity using the doubling of variable technique, see \cite[Remark p.18]{CL84}. In the absence of further assumptions, this rate is sharp in view of the recent examples discussed in \cite[Section 4]{Sprekeleretal}, and it becomes even slower if one weakens the regularity conditions on $H$ and $u$ or considers homogenization problems, see for example \cite{BCD,CCM}. 

One may then wonder whether the $\mathcal{O}(\sqrt{\eps})$ rate can be improved under further assumptions. First, following the common principle in PDE theory that smallness conditions imply smoothness, one can deduce a faster $\mathcal{O}(\eps)$ rate of convergence in sup-norms when $\|D^2u_T\|_{\infty}$ is \textit{small} or for \textit{short time horizons} \cite[Chapter 12]{L82Book}. This rate is a consequence of $C^{1,1}$ estimates that are independent of the viscosity, and that hold typically in these smallness regimes. In this case, applying the usual maximum principle to the (linearized) equation solved by $u_\eps-u$, one obtains
\[
\|u_\eps-u\|_{L^\infty(Q_T)}\leq C\eps.
\]
This agrees with $L^1$ rates of convergence for the vanishing viscosity process of the one-dimensional Burger's equation \cite[Theorem 3.1]{Wang}.\\

If $H$ is assumed to be \textit{uniformly convex}, well-known results provide \textit{one-side} $\mathcal{O}(\eps)$ rate in sup-norm \cite{L82Book} (this rate is also well-known for conservation laws under the Oleinik one-side Lipschitz condition \cite{NT}), and $\mathcal{O}(\eps)$ rate in the weaker $L^1$ norm \cite{LinTadmor,CGM}. These rates can be interpolated, giving new convergence results in $L^p$ spaces \cite{CGM}. A special role in the uniformly convex framework is played by the purely \textit{quadratic} Hamiltonian: by means of the Hopf-Cole reduction and the explicit representation formula for the heat flow \cite{Brenier}, the $\mathcal{O}(\eps |\log \eps|)$ rate in sup-norm has been observed in the literature. See \cite{Wang,TT} for conservation laws, and \cite[Proposition 4.4]{Sprekeleretal} for the one-dimensional Hamilton-Jacobi equation. This result has been very recently generalized to the multidimensional case \cite[Section 3]{ChaintronDaudin}, for equations posed on the whole space $\R^n$ and Lipschitz data. \\

The first goal of this paper is to prove that the $\mathcal{O}(\eps |\log \eps|)$ rate holds beyond the purely quadratic case, that is, for uniformly convex Hamiltonians. We assume, in particular, that for some $0<\theta\leq\Theta$, 
\begin{equation}\label{H}\tag{H}
\text{ $H\in C^2$ \, and \, $\theta \mathbb{I}_n\leq D^2_{pp}H(p)\leq \Theta \mathbb{I}_n$ for all $p$.}
\end{equation}

Throughout the paper, $u_T$ will be implicitly assumed to be of class $C^{2,\alpha}$, and $f$ of class $C^{2+\alpha, 1+\alpha/2}$. This will guarantee that the unique solution to \eqref{HJintroVisc} is classical for each $\eps > 0$. Nevertheless, all the estimates obtained below will depend on Lipschitz or semiconcavity properties of $u_T$ and $f$, that will be specified below. 
 
 The first result, stated in Theorem \ref{main1}, shows that the general $\mathcal{O}(\sqrt\eps)$ rate can be boosted to
\begin{equation}\label{power}
-C_{1,\beta}\eps^\beta\leq u_\eps-u\leq C_{2}\eps,\ \text{on $Q_T$, for any $\beta\in(1/2,1)$}.
\end{equation}
Here, $C_{1,\beta}, C_{2}$ depend on some semiconcavity properties of $u_T$ and $f$. This estimate implies also a new $\mathcal{O}(\eps^{\beta/2})$ rate of convergence of gradients in $L^\infty_t(L^2_x)$ norms and, consequently, the convergence rate for solutions of certain hyperbolic systems of conservation laws.

Secondly we prove, in Theorem \ref{main2}, the analytic bound
\begin{equation}\label{log}
\|(u_\eps-u)(\tau)\|_{L^\infty(\T^n)}\leq C_\tau \eps|\log\eps|, \quad \text{for any $\tau \in [0,T)$}.
\end{equation}
Note that $C_\tau$ deteriorates as $\tau \to T$, as we are using only the information $u_T\in W^{1,\infty}(\T^n)$ on the final condition. The bounds \eqref{power}-\eqref{log} complete the picture of the vanishing viscosity approximation for \textit{uniformly convex} Hamiltonians, answering a question posed in \cite[p.17]{Sprekeleretal}, which provided numerical evidence in favor of the validity of these results. Moreover, we observe in Remark \ref{linearrate} that the order $\eps$ can be reached in special situations, such as the one of convex data. \\

In the second part of the paper, we continue our analysis and consider Hamiltonians that are merely \textit{strictly convex}. More precisely, for $H\in C^2$ we relax \eqref{H} to 
\begin{equation}\label{H2}\tag{$H_\gamma$}
\text{$\mathrm{Tr}(D^2_{pp}H(p)M^2)\geq \theta |M|^2 |p|^{\gamma-2}$, $M\in \mathrm{Sym}_n,\ \gamma>2,\ $\text{ and } $D^2_{pp}H\leq\Theta I_n$, $|p|\leq C$,}
\end{equation}
for some $0<\theta <\Theta$, a model being the superquadratic Hamiltonian
\[
H(p)=|p|^\gamma,\ \gamma>2.
\]
We prove in Theorem \ref{main3} the following estimate: there exists $\beta_{\gamma}\in\left(\frac12,1\right)$ such that
\begin{equation}\label{powergamma}
- C_{3,\beta}\,\eps^{\beta}\leq u_\eps-u\leq C_{2}\, \eps \quad \text{on $Q_T$, for each $\beta < \beta_\gamma$} .
\end{equation}
The upper bound $\beta_{\gamma}$ is explicit. One can check that the leftmost rate $\mathcal{O}(\eps^{\beta_{\gamma}})$ deteriorates to $\mathcal{O}(\sqrt{\eps})$ as $\gamma\to+\infty$ and agrees with \eqref{power} in the limit $\gamma\to2$. To the best of our knowledge, these intermediate rates $\beta_\gamma$ are new within the framework of strictly convex Hamiltonians. Numerical experiments have been carried out in \cite{Sprekeleretal}, suggesting that one should indeed expect some intermediate situation $\beta_\gamma \in (1/2,1)$. We do not known whether the upper bound $\beta_\gamma$ obtained here is sharp or not, and we believe that further investigation of strictly but not uniformly convex problems can be an interesting research direction. See Remark \ref{superquadrem} for further comments.\\

The proofs  presented here exploit a duality approach \cite{EvansARMA}, being based on the analysis of the adjoint of the linearization of \eqref{HJintroVisc}. A more detailed description of the arguments is given in the next section. 

We conclude by mentioning that one of our main motivations to address this vanishing viscosity problem was the more complicated study of the rate of the convergence problem in Mean Field Control, where, roughly speaking, one fixes $\eps > 0$ in \eqref{HJintroVisc} and let the dimension $n \to \infty$. In some specific situations (see for example \cite{DaudinDelarueJackson}), this problem can be reduced in fact to a question of vanishing viscosity for Hamilton-Jacobi equations posed on the Euclidean space (with fixed dimension). We report the implications of our results in the convergence problem for Mean Field Control in Remark \ref{mfc}.

It is worth mentioning that other approximations of Hamilton-Jacobi equations have been considered in the literature, for example regularizations via nonlocal operators \cite{DroniouImbert, Goffi}. Finally, it would be interesting to apply this approach to the study of the convergence rate of the vanishing viscosity approximation of Mean Field Games systems \cite{TangZhang}, numerical methods, cf. \cite{BCD,Calder,CGT}, stationary problems and homogenization. 
\bigskip

\textbf{Acknowledgements.} The authors would like to thank H. V. Tran and F. Meng for some comments on the first draft of the manuscript and the anonymous referee for the careful reading. The authors are member of the Gruppo Nazionale per l'Analisi Matematica, la Probabilit\`a e le loro Applicazioni (GNAMPA) of the Istituto Nazionale di Alta Matematica (INdAM), and they are partially supported by the INdAM-GNAMPA projects 2025. M.C. has been partially funded by King Abdullah University of Science and Technology Research Funding (KRF) under award no. ORA-2021-CRG10-4674.2,  and by EuropeanUnion–NextGenerationEU under the National Recovery and Resilience Plan (NRRP), Mission 4 Component 2 Investment 1.1 - Call PRIN 2022 No. 104 of February 2, 2022 of Italian Ministry of University and Research; Project 2022W58BJ5 (subject area: PE - Physical Sciences and Engineering) “PDEs and optimal control methods in mean field games, population dynamics and multi-agent
models”.

\section{Preliminaries and method of proof}

Following a classical approach (see for instance \cite[Theorem 5.1]{FlemingSouganidis}) we first observe that the estimate of the quantity $\|u_\eps-u\|_{L^\infty}$ can be obtained, via an asymptotic expansion, from an estimate of $\|\partial_\eps u_\eps\|_{L^\infty}$, see in particular Remark \ref{concl} below. 
Note also that, by differentiating \eqref{HJintroVisc} with respect to $\eps$, $v^\eps = \partial_\eps u_\eps$  solves the linear drift-diffusion PDE
\[
-\partial_t v^\eps-\eps\Delta v^\eps+D_pH(Du_\eps)\cdot Dv^\eps=\Delta u_\eps 
\]
equipped with homogenous final condition $v^\eps(T)=0$. Now, for any $(\bar x,\tau) \in Q_T$, by means of the dual problem solved by $\rho = \rho_{\eps, \tau, \bar x}$
\begin{equation}\label{fpintro}
\begin{cases}
\partial_t \rho-\eps\Delta \rho-\mathrm{div}(D_pH(Du_\eps)\rho)=0&\text{ in }Q_{\tau,T}:=\T^n\times(\tau,T),\\
\rho(\tau)=\delta_{\bar x}&\text{ in }\T^n,
\end{cases}
\end{equation}
we get a representation formula, for arbitrary $s \in [\tau, T]$, as follows
\[
v^\eps(\bar x, \tau)=\int_{\T^n}v^\eps(s)\rho(s)\,\dd x+\iint_{\T^n\times(\tau,s)}\Delta u_\eps\rho\dd x \dd t=\mathrm{(I)}+\mathrm{(II)} \, .
\]
One can take for simplicity $s = T$, so that the integral (I) is zero. Therefore, the estimate on $v^\eps$ boils down to a control on (II), that involves second-order derivatives $D^2 u_\eps$ \textit{weighted} by $\rho$. The crucial step of the approach presented here is in fact to exploit the weight $\rho$, which, from perspective of optimal control, describes the flow of optimal trajectories originating from $\bar x$ at time $\tau$.

There are two main ingredients that allow to control the quantity $\iint \Delta u_\eps \rho$. These are the presence of  the diffusion and the convexity of the Hamiltonian. The diffusion gives bounds of the form
\begin{equation}\label{diffest}
\int_\tau^T \int_{\T^n} |D^2u_\eps|^2\rho\dd x \dd t\lesssim \frac{1}{\eps},
\end{equation}
see Lemma \ref{lipbound} below. Note that this implies immediately, by H\"older's inequality and conservation of mass for $\rho$, that
\[
\left| \iint_{\T^n\times(\tau,T)}\Delta u_\eps\rho \right| \le \sqrt{\frac{T-\tau}{\eps}}.
\]
This crude bound gives the $\mathcal O(\sqrt{\eps})$ control on $\|u_\eps-u\|_{L^\infty}$. To improve it, one needs to get further information out of the convex structure of $H$. To this aim, we follow a general principle of optimal control, which roughly states that optimal trajectories should be ``stable'', as well as the optimal drift, when restricted to time intervals of the form $[t, T] \subset (\tau, T]$, that is, away from the originating point $(\bar x, \tau)$. We make this observation quantitative by proving estimates of the form
\begin{equation}\label{convexest}
\int_{\tau+\delta}^T \int |D^2u_\eps|^2\rho\dd x \dd t\lesssim \frac{1}{\delta^\alpha},
\end{equation}
where $\delta > 0$ and $\alpha > 1$ (see \eqref{secondord}). This, together with \eqref{diffest}, allows to improve the previous bounds, by splitting the time interval into $[\tau, \tau+ \delta] \cup [\tau+ \delta, T]$ as follows:
\[
\left| \iint_{\T^n\times(\tau,T)}\Delta u_\eps\rho \right| \lesssim \frac{1}{\delta^\alpha} + \sqrt{\frac{\delta}{\eps}}.
\]
A suitable choice of $\delta = \delta(\eps)$ gives then the desired rate $\mathcal O(\eps^\beta)$, for $\beta < 1$.

To reach the optimal rate $\mathcal O(\eps|\log\eps|)$, a more delicate analysis is needed, as it corresponds to the ``critical'' case $\alpha = 1$ in the previous bound. We exploit here some additional entropy and first-order bounds for \eqref{fpintro} when the initial datum does not belong to the Orlicz space $L {\log}L$ (Lemma \ref{lemDrho}). We refer to \cite{NSS} for a related analysis under different assumptions on the velocity field and the initial datum. \\

The analysis of the strictly convex case is even more delicate, as \eqref{convexest} is meaningful only when $|Du_\eps|$ is bounded away from zero. Under our strict convexity assumptions, \eqref{convexest} becomes
\[
\int_{\tau+\delta}^T \int_{|Du|^2 \ge m} |D^2u_\eps|^2\rho\dd x \dd t\lesssim \frac{1}{m^{\frac{\gamma-2}2}\delta^\alpha},
\]
To compensate the lack of information when $|Du|$ is close to zero, we need to sharpen \eqref{diffest} in regions where the gradient is small. The proof of the crucial estimate \eqref{so2} is obtained via a Bochner-type identity for $\tilde w=\varphi(|Du|^2)$,
\[
\Delta \tilde w=4\varphi''\sum_{j=1}^n(Du_{x_j}\cdot Du)^2+2\varphi'D\Delta u\cdot Du+2\varphi'|D^2u|^2,
\]
where $\varphi$ is a carefully chosen truncation function. We proceed with a hole-filling technique in the spirit of K.-O. Widman \cite{Widman}. This technique has been historically applied to obtain interior Morrey estimates, and local H\"older continuity as a byproduct, for elliptic equations and systems. We apply here a similar method to get integral bounds of (weighted) second order derivatives over sublevel sets of the gradient (instead of concentric balls as in the classical literature).

We conclude by stressing again that it is crucial to look at integral estimates on second derivatives, as classical pointwise estimates (coming from example from the theory of viscosity solutions) merely give
\[
\|\partial_t u_\eps\|_{L^\infty_{x,t}},\|Du_\eps\|_{L^\infty_{x,t}}\leq C\implies \|\Delta u_\eps\|_{L^\infty_{x,t}}\lesssim \frac{1}{\eps}.
\]
These bounds have no (straightforward) use in the control of the rate of convergence. \\

The next lemmata make some of the steps that have been previously described rigorous.

\begin{lemma}\label{eqdereps} Assume that $H$ satisfies \eqref{H} or \eqref{H2}. For every $(x,t) \in \T^n\times[0,T]$ and $\eps > 0$, the derivative $v^\eps(x,t) = \partial_\eps u_\eps(x,t)$ is well-defined. The function $v^\eps$ is of class $C^{2+\alpha,1+\alpha/2}(\T^n\times[0,T])$, and it solves 
\begin{equation}\label{eqv}
\begin{cases}
-\partial_t v^\eps-\eps\Delta v^\eps+D_pH(Du_\eps)\cdot Dv^\eps=\Delta u_\eps &\text{ in }\T^n\times(0,T),\\
v^\eps(x,T)=0&\text{ in }\T^n
\end{cases}
\end{equation}
in the classical sense.
\end{lemma}

\begin{proof} We denote by $C^{2k+\alpha,k+\alpha/2}$ the classical H\"older space with respect to the parabolic distance. Note first that by Schauder estimates (for quasilinear parabolic equations), $u_\eps$ is of class $C^{2+\alpha,1+\alpha/2}(\T^n\times[0,T])$, hence $\Delta u_\eps \in C^{\alpha,\alpha/2}(\T^n\times[0,T])$ and the (linear) Cauchy problem \eqref{eqv} admits a unique classical solution $v^\eps \in C^{2+\alpha,1+\alpha/2}(\T^n\times[0,T])$. Consider now, for $\eta > 0$, the finite difference
\[
v^\eps_\eta=\frac{u_{\eps+\eta}-u_\eps}{\eta}
\]
and note that it satisfies $v^\eps_\eta(T) = 0$ and
\[
-\partial_t v^\eps_\eta-(\eps+\eta)\Delta v^\eps_\eta+\frac{H(Du_{\eps+\eta})-H(Du_\eps)}{\eta}=\Delta u_\eps.
\]
Here we use the regularity of $H$ to show that, by the second-order Taylor expansion,
\[
H(Du_{\eps+\eta})=H(Du_\eps)+D_pH(Du_\eps)\cdot D(u_{\eps+\eta}-u_\eps)+  \Psi_{\eps,\eta} D(u_{\eps+\eta}-u_\eps)\cdot D(u_{\eps+\eta}-u_\eps),
\]
where the matrix-valued function $\Psi_{\eps,\eta} $ is defined as
\[
\Psi_{\eps,\eta}(x,t) = \frac 12 \int_0^1 (1-\zeta) D^2_{pp} H\big(Du_\eps(x,t) + \zeta(Du_{\eps+\eta}(x,t)-Du_{\eps}(x,t))\big) \dd \zeta.
\]
The function $\Psi_{\eps,\eta}$ is continuous, and bounded in $\T^n \times [0,T]$ uniformly with respect to $\eta$ by the uniform bound $D^2_{pp}H\leq \Theta$. Hence,
\[
-\partial_t v^\eps_\eta-(\eps+\eta)\Delta v^\eps_\eta+D_pH(Du_\eps)\cdot Dv^{\eps}_\eta + \eta \Psi_{\eps,\eta} Dv^\eps_\eta \cdot Dv^\eps_\eta =\Delta u_\eps.
\] 

To pass to the limit $\eta \to 0$, we use the H\"older regularity estimates on $v^\eps_\eta$ from \cite{LSU} (uniform in $\eta$), that go through gradient bounds first, then $W^{2,1}_p$ (the parabolic Sobolev space of second order) regularity and Sobolev embeddings finally. We can then extract a subsequence, that we keep denoting by $v^\eps_\eta$, that converges uniformly to some continuous function $v^\eps$ by the Ascoli-Arzel\`a Theorem. To show that the limit $v^\eps$ (along the specific subsequence) solves \eqref{eqv}, we use a standard viscosity solution argument: $\eta \Psi_{\eps,\eta}$ vanishes uniformly on $\T^n \times [0,T]$ and $\eps + \eta \to \eps$. Therefore, by the stability of viscosity solutions, $v^\eps$ is a viscosity solution of \eqref{eqv}. Since viscosity solutions to the linear problem  \eqref{eqv} are unique, $v^\eps$ enjoys in fact the regularity claimed at the beginning of the proof. Moreover, this uniqueness property also shows that $v^\eps_\eta \to v^\eps$ uniformly on $\T^n \times [0,T]$ along the full limit $\eta \to 0$, providing the desired assertion.
\end{proof}

\begin{rem}\label{concl}
The proof of the rate of convergence will be reduced to an $L^\infty$ estimate on $v^\eps$. Having this at our disposal, we then get that, for every $\eps'' > \eps' > 0$, $(x,t) \in \T^n \times [0,T]$,
\[
|u_{\eps''}(x,t) - u_{\eps'}(x,t)| \le \int_{\eps'}^{\eps''} |v^\eps(x,t)| \dd \eps \le \int_{\eps'}^{\eps''} \|v^\eps\|_\infty \dd \eps.
\]
Hence, whenever $\eps \mapsto \|v^\eps\|_\infty$ is integrable on $(0,1)$, we have that $u_\eps$ is a Cauchy sequence in the space of continuous functions on $\T^n \times [0,T]$, and its uniform limit $u$ satisfies
\[
\|u_{\eps''} - u\|_\infty \le \int_{0}^{\eps''} \|v^\eps\|_\infty \dd \eps.
\]
Moreover, by uniform convergence and stability of the notion of viscosity solution, $u$ turns out to be a (the) viscosity solution of \eqref{HJintro}. Similarly, an upper bound on $(v^\eps)^+$ yields an upper bound on $(u_{\eps} - u)^+$.
\end{rem}

\begin{lemma}\label{lipbound} Assume that $H$ satisfies \eqref{H} or \eqref{H2} (in fact, $H$ can be assumed to be only locally Lipschitz). Then, for each $(x, \tau) \in \T^n \times [0,T]$ and $\rho$ solving \eqref{fpintro} we have
\[
|Du_\eps(x, \tau)| + \eps \int_\tau^T\int_{\T^n}|D^2u_\eps|^2\rho\dd x \dd t\leq C_L,
\]
where $C_L$ depends on $\|D u_T\|_{L^\infty_x} , \|D f\|_{L^{\infty}_{x,t}}$ and $T$.
\end{lemma}

\begin{proof} This estimate, first observed by L.C. Evans \cite{EvansARMA}, is well-known, but we briefly sketch its proof here for completeness. It is enough to look at the equation solved by $w=|Du_\eps|^2$, namely
\[
-\partial_t w-\eps\Delta w+2\eps|D^2u_\eps|^2+DH (Du_\eps)\cdot Dw=Df\cdot Du_\eps.
\] 
By duality with \eqref{fpintro}, we have
\begin{multline*}
|Du_\eps(x, \tau)|^2 + 2\eps \int_\tau^T\int_{\T^n}|D^2u_\eps|^2\rho\dd x \dd t \le \int_{\T^n}|D u_T|^2\rho(x,0)\dd x + \int_\tau^T\int_{\T^n}Df\cdot Du_\eps \rho\dd x \dd t \\
\le \|D u_T\|^2_{L^\infty(\T^n)} +  T \|D f\|_{L^\infty(Q_T)}\|D u\|_{L^\infty(Q_T)}.
\end{multline*}
Using now Young's inequality and taking the supremum over $(x, \tau) \in \T^n \times [0,T]$ one concludes.
\end{proof}

\section{Rate of convergence for uniformly convex Hamiltonians and semiconcave data}\label{sec;power}
In this section we prove that the rate of convergence of the quantity $\|(u_\eps-u)^-\|_{L^\infty}$ for semiconcave solutions of Hamilton-Jacobi equation can be boosted from $\mathcal{O}(\sqrt{\eps})$ to $\mathcal{O}(\eps^\beta)$, $\beta\in(1/2,1)$ in $L^\infty$ norm. This is usually known to be $\mathcal{O}(\sqrt\eps)$ in $L^\infty$ norm or $\mathcal{O}(\eps)$ in the weaker $L^1$ norm or in a certain average sense \cite{TranBook}, see also \cite{Goffi} for different rates of order $\mathcal{O}(\eps)$ with nonlocal regularization of the PDE. The proof uses the nonlinear adjoint method and exploits only the conservation of mass property of the dual equation \eqref{fpintro}.

\begin{thm}\label{main1}
Let $H$ be such that \eqref{H} holds. 
Then, if $u_\eps$ solves \eqref{HJintroVisc} and $u$ solves the first-order equation \eqref{HJintro} we have that for all $\beta\in\left(\frac12,1\right)$ there exists $C=C(\beta)>0$, which depends also on $n,\theta,\|(\Delta u_T)^+\|_{L^\infty_x}$, $\|(\Delta f)^+\|_{L^\infty_{x,t}}, T, \|Du_T\|_{L^\infty_x},\|Df\|_{L^\infty_{x,t}}$ such that
\[
-C\eps^\beta\leq u_\eps-u\leq (T\|(\Delta u_T)^+\|_{L^\infty_x}+\|(\Delta f)^+\|_{L^1_t(L^\infty_x)})\eps,\quad \beta\in(1/2,1).
\]
\end{thm}

\begin{proof}
\textit{\underline{Step 1} ($\mathcal{O}(\eps)$ bound from above)}. The linear upper bound on $(u_\eps-u)^+$ is well-known, even under weaker assumptions on $H$, see \cite{CGM} and \cite[Proposition 11.2]{L82Book}. We briefly recall the proof for completeness. First, note that by Lemma \ref{eqdereps} the function $v^\eps=\partial_\eps u_\eps$ solves \eqref{eqv}. Testing \eqref{eqv} against the solution of the adjoint problem \eqref{fpintro} , we find by the conservation of mass for \eqref{fpintro} that
\begin{multline*}
v^\eps(\bar x, \tau) =\int_\tau^T\int_{\T^n}\Delta u_\eps\rho\dd x \dd t\leq \|(\Delta u_\eps)^+\|_{L^1_t(L^\infty_x)}\\
\leq T\|(\Delta u_T)^+\|_{L^\infty(\T^n)}+\|(\Delta f)^+\|_{L^1(0,T;L^\infty(\T^n))},
\end{multline*}
where we used the semi-superharmonic bound in \cite[Remarks 3.6 and 4.9]{CGM}, see also \cite{L82Book}. By Remark \ref{concl} we conclude the estimate
\[
\|(u_\eps-u)^+\|_{L^\infty(Q_T)}\leq \eps(T\|(\Delta u_T)^+\|_{L^\infty(\T^n)}+\|(\Delta f)^+\|_{L^1(0,T;L^\infty(\T^n))}).
\]

\medskip
\textit{\underline{Step 2} (Second order estimate)}. We prove that for $\alpha\in(1,2)$ we have the bound
\begin{equation}\label{secondord}
\int_\tau^T\int_{\T^n}(t-\tau)^\alpha|D^2u_\eps|^2\rho\dd x \dd t\leq K,
\end{equation}
where $\rho$ solves \eqref{fpintro} and $K$ depends on $n, \alpha, T,\theta,u_T,f$. To this aim, we exploit the uniform convexity of $H$, and follow an argument similar to that described in  \cite{CGMT,LMT}. We first find, by differentiating twice the equation for $u_\eps$ and the uniform convexity of $H$, the following inequality solved by the function $z(x,t)=z_\eps(x,t)=(t-\tau)^\alpha  \Delta u_\eps(x,t)$:
\begin{multline*}
-\partial_t z-\eps\Delta z+\theta (t-\tau)^\alpha  |D^2u_\eps|^2 +D_pH(Du_\eps)\cdot Dz\\
\leq -\alpha (t-\tau)^{\alpha-1}\Delta u_\eps(x,t)+(t-\tau)^\alpha\Delta f(x,t) \quad \text{ in }Q_T.
\end{multline*}
By duality and integrating in $\T^n\times(\tau,T)$ we have
\begin{multline}\label{secondorde}
\underbrace{\int_{\T^n}z(\tau)\rho(\tau)\,dx}_{=0}+\theta\int_\tau^T\int_{\T^n}(t-\tau)^\alpha |D^2u_\eps|^2 \rho\dd x \dd t= \underbrace{\int_{\T^n}z(T)\rho(T)\dd x}_{\leq (T-\tau)^\alpha\|(\Delta u_T)^+\|_{L^\infty(\T^n)}}\\
-\alpha \int_\tau^T\int_{\T^n} (t-\tau)^{\alpha-1}\Delta u_\eps \rho\dd x \dd t+\int_\tau^T\int_{\T^n} (t-\tau)^{\alpha} \Delta f \rho\dd x \dd t.
\end{multline}
We now use Young's inequality as follows
\begin{align*}
-\alpha &\int_\tau^T\int_{\T^n} (t-\tau)^{\alpha-1} \Delta u_\eps \rho\dd x \dd t\\
&\leq \frac{\theta}{2n}\int_\tau^T\int_{\T^n}(t-\tau)^\alpha|\Delta u_\eps |^2\rho\dd x \dd t+\frac{n\alpha^2}{2\theta}\int_\tau^T\int_{\T^n}(t-\tau)^{\alpha-2}\rho\dd x \dd t\\
&\leq \frac{\theta}{2}\int_\tau^T\int_{\T^n}(t-\tau)^\alpha|D^2u_\eps|^2\rho\dd x \dd t+\frac{n\alpha^2}{2\theta}\int_\tau^T (t-\tau)^{\alpha-2}\,dt\\
&\leq\frac{\theta}{2}\int_\tau^T\int_{\T^n}(t-\tau)^\alpha|D^2u_\eps|^2\rho\dd x \dd t+\frac{n\alpha^2}{2\theta(\alpha-1)} T^{\alpha-1}.
\end{align*}
We then obtain
\begin{multline}\label{est1}
\int_\tau^T\int_{\T^n}(t-\tau)^\alpha|D^2u_\eps|^2\rho\dd x \dd t\\
\leq \frac{n \alpha^2}{\theta^2(\alpha-1)}T^{\alpha-1}+\frac{2}{\theta}T^\alpha\|(\Delta u_T)^+\|_{L^\infty(\T^n)}+\frac{T^{\alpha+1}}{\alpha+1}\|(\Delta f)^+\|_{L^\infty(Q_T)}=:K.
\end{multline}
\medskip
\textit{\underline{Step 3} ($\mathcal{O}(\eps^\beta)$ bound from below).} We now proceed to estimate the rate of convergence: define again $v^\eps(x,t) = \partial_\eps u_\eps(x,t)$ and find by Lemma \ref{eqdereps} the PDE \eqref{eqv}. By duality we find
\begin{equation}\label{repr}
v^\eps(\bar x,\tau)=\int_{\T^n}v^\eps(T)\rho(T)\dd x+\int_\tau^T\int_{\T^n}\Delta u_\eps\rho\dd x \dd t.
\end{equation}
We only estimate the last term on the right-hand side since $v^\eps(T)=0$. Assume first that $\tau + \eps < T$. Then,
\begin{align*}
\Big| \int_\tau^T\int_{\T^n}\Delta u_\eps\rho\dd x \dd t \Big| &\leq \sqrt{n}\left(\int_{\tau + \eps}^{T}\int_{\T^n}(t-\tau)^{\alpha/2}|D^2u_\eps|\rho (t-\tau)^{-\alpha/2}\dd x \dd t+\int_{\tau}^{\tau + \eps}\int_{\T^n}|D^2u_\eps|\rho\dd x \dd t\right)\\
&\le \sqrt{n}\left(\int_{\tau + \eps}^{T}\int_{\T^n}(t-\tau)^{\alpha}|D^2 u_\eps|^2\rho\dd x \dd t\right)^\frac12\left(\int_{\tau + \eps}^{T}\int_{\T^n} (t-\tau)^{-\alpha}\rho\dd x \dd t\right)^\frac12\\
& \quad + \sqrt{n} \left(\int_{\tau}^{\tau + \eps}\int_{\T^n}|D^2u_\eps|^2\rho\dd x \dd t \right)^\frac12\left(\int_{\tau}^{\tau + \eps}\int_{\T^n}\rho\dd x \dd t \right)^\frac12\\
&\leq \sqrt{\frac{nK\eps^{1-\alpha}}{\alpha-1}} + \sqrt {nC_L},
\end{align*}
where we used the estimate \eqref{est1} and Lemma \ref{lipbound}. Note that if $\tau \ge T-\eps$, then the same estimate holds, without the term $\sqrt{\frac{nK\eps^{1-\alpha}}{\alpha-1}}$, as in this case there is no need to split the first integral.

Back to \eqref{repr}, this implies that 
\[
v^\eps(\bar x, \tau) \ge - \sqrt{\frac{nK}{\alpha-1}} \eps^{\frac{1-\alpha}2} - \sqrt {nC_L}.
\]
Since the previous estimate holds for all $\bar x, \tau$, by Remark \ref{concl} we obtain, by setting $\frac{3-\alpha}2=\beta\in(1/2,1)$, 
\[
u_{\eps}-u\geq - \frac1\beta\sqrt{\frac{nK}{2(1-\beta)}} \eps^{\beta} - \sqrt {nC_L} \eps.
\]
\end{proof}

\begin{rem} Since $K$ in the previous proof grows linearly in the dimension $n$, one can easily check that the rate (from below) of the previous result grows linearly in $n$.

The previous proof extends to more general Hamiltonians $H=H(x,t,p)\in C^2(\R^n\times(0,T)\times\R^n)$ under the conditions
\[
\sup_{(x,t)\in Q_T}|D_x H(x,t,p)|, \sup_{(x,t)\in Q_T}|D^2_{xp} H(x,t,p)|, \sup_{(x,t)\in Q_T}|D^2_{xx} H(x,t,p)|\leq C_H(1+|p|),
\]
together with
\[
D^2_{pp}H(x,t,p)\xi\cdot \xi\geq a(x,t)|\xi|^2,\ a>0,
\]
and the constant of the estimates will depend also on $C_H,a$. More precisely, the bound \eqref{secondord} under these assumptions was proved in \cite{TranCalcVarPDE}, while the second order estimate can be found starting from equation (36) in \cite{CGsima} and following the proof of Proposition 3.6 therein.
\end{rem}

\begin{rem}[Reaching $\mathcal O(\eps)$]\label{linearrate} As we mentioned in the introduction, the linear order $\mathcal{O}(\eps)$ can be reached in exceptional cases only, for example under smallness conditions. We mention here two other possible situations where the presence of a uniformly convex Hamiltonian can be exploited. If one takes $\alpha = 0$ in the previous proof, then \eqref{secondorde} reads
\[
\theta\int_\tau^T\int_{\T^n} |D^2u_\eps|^2 \rho\dd x \dd t \le \|(\Delta u_T)^+\|_{L^\infty(\T^n)} +
\int_\tau^T\int_{\T^n}\Delta f \rho\dd x \dd t - \Delta u_\eps(\bar x, \tau).
\]
The last term $- \Delta u_\eps(\bar x, \tau)$ is in general unbounded as $\eps \to 0$, and that is exactly why we put the weight $(T-\tau)^\alpha$ in the previous proof. Nevertheless, if one assumes that $u_T$ and $f(\cdot, t)$ are \textit{convex}, then it is known that $u_\eps(\cdot, t)$ is convex as well for all $t$ and $\eps$ (for example by employing coupling methods, see for example \cite{conforti} or \cite{alpar}). Therefore, $- \Delta u_\eps(\bar x, \tau) \le 0$ and the previous inequality gives a uniform bound on $\iint |D^2u_\eps|^2 \rho$ that is independent of $\tau, \bar x, \eps$. Arguing as in the previous proof (with an even simpler argument, since there is no need to split the time interval), one reaches the order $\mathcal{O}(\eps)$ in sup-norm. Clearly, the hypothesis on convexity of the data is extremely restrictive if periodicity is also assumed, but the arguments presented here can be easily adapted to the nonperiodic setting, paying some attention on the usual technicalities arising from the presence of unbounded domains.

In fact, even if the data are not assumed to be convex, every point $(\bar x, \tau)$ such that $- \Delta u_\eps(\bar x, \tau)$ remains bounded above uniformly with respect to $\eps$ is a point where convergence of $u_\eps(\bar x, \tau)$ to $u(\bar x, \tau)$ is linear in $\eps$. The points where convergence is linear are known to form a set of full Lebesgue measure when the Hamiltonian is purely quadratic, see \cite[Proposition 4-(ii)]{Sprekeleretal}. We expect also this to be true when $H$ is uniformly convex.
\end{rem}

The first consequence of the quantitative bound in Theorem \ref{main1} concerns the speed of convergence of gradients in the vanishing viscosity approximation.

\begin{cor}\label{convgrad}
Under the assumptions of Theorem \ref{main1} we have that for all $\beta\in\left(\frac12,1\right)$ there exists $C=C(\beta)>0$, which depends also on $n,\theta,\|(\Delta u_T)^+\|_{L^\infty_x}$,$\|(\Delta f)^+\|_{L^\infty_{x,t}}, T, \|Du_T\|_{L^\infty_x},\|Df\|_{L^\infty_{x,t}}$ such that
\[
\|Du_\eps-Du\|_{L^\infty(0,T;L^2(\T^n))}\leq C\eps^\frac{\beta}{2}.
\]
\end{cor}
\begin{proof}
We first note that the semiconcavity estimate implies
\[
\|\Delta u_\eps\|_{L^\infty(0,T;L^1(\T^n))}\leq C_1
\]
for a constant $C_1$ independent of $\eps$, see \cite{K66d} or \cite[Theorem 4.14]{CGM}. Integrating by parts we get for a.e. $t\in(0,T)$
\begin{multline*}
\int_{\T^n}|Du_\eps(t)-Du(t)|^2\dd x=-\int_{\T^n}\Delta(u_\eps-u)(u_\eps-u)(t)\dd x\\
\leq \|(u_\eps-u)(t)\|_{L^\infty(\T^n)}\|\Delta(u_\eps-u)(t)\|_{L^1(\T^n)}\leq C\eps^\beta,
\end{multline*}
where $C$ is independent of $\eps$.
\end{proof}

We conclude this section with an application to the convergence rate of the solution to the quasilinear parabolic system
\[
\begin{cases}
\partial_t p_i^\eps+\partial_{x_i}H(p_1^\eps,...,p_n^\eps)=\eps\Delta p_i^\eps&\text{ in }Q_T,\ i=1,...,n\\
p_i^\eps(x,0)=p_i^0(x)&\text{ in }\T^n
\end{cases}
\]
to the solution of the hyperbolic system
\begin{equation}\label{hyperbolic}
\begin{cases}
\partial_t p_i+\partial_{x_i}H(p_1,...,p_n)=0&\text{ in }Q_T,\ i=1,...,n\\
p_i(x,0)=p_i^0(x)&\text{ in }\T^n.
\end{cases}
\end{equation}
The convergence of the parabolic equation to the system of conservation laws \eqref{hyperbolic} was proved by S.N. Kruzhkov in Theorem 8 of \cite{K67II}, see also \cite[Section 16.1]{L82Book}. The connection among \eqref{HJintro} and \eqref{hyperbolic} is the following: if $u$ is a solution of \eqref{HJintro}, then $p(x,t)=Du(x,t)$ with $p_i(x,t)=\partial_{x_i}u(x,t)$ and $p(x,0)=Du(x,0)$ lead to \eqref{hyperbolic}. Corollary \ref{convgrad} implies the following:
\begin{cor}
If $H$ is uniformly convex and $p_0=Du_0$ with $p_0\in C^1(\T^n)$, we have that for all $\beta\in\left(\frac12,1\right)$ there exists $C=C(\beta)>0$, which depends also on $n,\theta,\|(\mathrm{div}(p_0))^+\|_{L^\infty_x}, T, \|p_0\|_{L^\infty}$ such that
\[
\|p_i^\eps-p_i\|_{L^2(Q_T)}\leq C\eps^{\frac{\beta}{2}},\ \beta\in(1/2,1).
\]
\end{cor}

\section{An improved rate for uniformly convex Hamiltonians and  Lipschitz data}\label{sec;log}
We now improve the left-side rate of convergence of Section \ref{sec;power} to $\mathcal O(\eps |\log \eps|)$. While in the previous result we rely on the semiconcavity of $u_T$, the following argument will depend merely on its Lipschitz continuity, at the price of an estimate that deteriorates close to $T$. This provides an analytic proof to general uniformly convex Hamilton-Jacobi equations of the one-dimensional result in Proposition 4.4-(i) of \cite{Sprekeleretal}. 

\begin{thm}\label{main2}
Let $u_\eps$ be a solution of \eqref{HJintroVisc} with $H$ satisfying \eqref{H}. For all $\eps \in (0,1)$ and  $\tau \in [0, T)$ we have
\[
\|(u_\eps-u)(\tau)\|_{L^\infty(\T^n)}\leq \frac{C}{(T-\tau)^{1/2}}\eps|\log\eps|
\]
where $C$ depends on $n, \|Du_T\|_{L^\infty_x}, T, \|(\Delta f)^+\|_{L^\infty_{x,t}}, \theta, \Theta$.
\end{thm}
The proof of this result is divided in several steps and requires a careful analysis of entropy-type bounds for solutions of the adjoint problem \eqref{fpintro}. The next lemma shows a certain control with respect to $\eps$ of the $W^{1,1}$ norm for solutions of Fokker-Planck equations with bounded velocity fields and Dirac initial condition, in the regime of small viscosity. These estimates are known when the drift of the Fokker-Planck equation belongs to some Sobolev space and the initial data belongs to an Orlicz class,  cf. \cite{LBL}. In particular, this complements the analysis of \cite{NSS}, see Remark 1 therein, where a decay result of the Sobolev norm is studied when the velocity field $b=b(x,t)$ of the advection equation is weakly compressible, i.e. 
\[
[\mathrm{div}(b)]^-\in L^1_t(L^\infty_x),
\]
and $\rho(\tau)\in \mathrm{LLogL}(\T^n)$. These available results cannot be applied in our setting since here $\rho(\tau)$ is just a probability measure. Quantitative convergence rates for Fokker-Planck equations within the Diperna-Lions framework have been the matter of \cite{BCC}, see also the references therein. 

\begin{lemma}\label{lemDrho} Assume that $0 <\eps \le 1$ and $\tau \le T-4\eps$. Let $b\in L^\infty(\T^n \times (\tau, T))$ and $\rho:=\rho_{\eps, \tau, x}$ be a solution of
\[
\begin{cases} 
\partial_t \rho - \varepsilon \Delta \rho + { \rm div}(b \rho) = 0 & \text{ in } Q_{\tau,T}:=\T^n \times (\tau, T), \\
\rho(\tau) = \delta_x &
\end{cases}
\]
and assume that $\|b\|_{L^\infty(\T^n \times (\tau, T))} \le K$. Then, there exists $C>0$ depending on $n,K$
and
\[
t_1 \in [\tau+\eps, \tau + 2\eps], \qquad t_2 \in [T-2\eps, T-\eps]
\]
 such that
 \[
 \int_{\T^n} |D \rho(y,t_i)| \dd y \le \frac{C}{\eps}\left(1 + |\log \eps|\right)^{1/2}, \qquad i=1,2.
 \]

\end{lemma}

\begin{proof} After the time change $t \mapsto t-\tau$ the problem reads
\[
\begin{cases}
\partial_t \hat{\rho} - \varepsilon \Delta \hat{\rho} + {\rm div}(b \hat{\rho}) = 0 & \text{in } \T^n \times (0, T-\tau), \\
\hat{\rho}(0) = \delta_x &
\end{cases}
\]
We can apply \cite[Corollary 7.2.3]{BKRS} with the choice $A = \varepsilon \mathbb{I}_n$, $B = \hat b$, $c = 0$, $\lambda_0 = \varepsilon$, $\Theta = \frac{1}{2}$, and any $\nu > (n+2)/2$ to get
\begin{multline*}
\hat \rho(x, t) \leq C_{\nu, n} \left( 1 + \frac{1}{\varepsilon} \right)^\nu t^{-\frac{n+2}{2}} 
\int_{t/2}^t \int_{\T^n} \left[ 1 + \varepsilon^\nu + \frac{t^{2\nu}}{\varepsilon^\nu} |\hat b|^{2\nu} \right] \hat \rho \, \dd y \, \dd s \\
\leq C_{\nu, n, K} \left( 1 + \frac{1}{\varepsilon} \right)^\nu t^{-\frac{n+2}{2}}   \left[ 1 + \frac{t^{2\nu}}{\varepsilon^\nu} \right] \int_{t/2}^t \int_{\T^n} \hat \rho \, \dd y \, \dd s
= C_{\nu, n, K} \left( 1 + \frac{1}{\varepsilon} \right)^\nu t^{-\frac{n}{2}}   \left[ 1 + \frac{t^{2\nu}}{\varepsilon^\nu} \right] 
\end{multline*}
for all $(x,t) \in \T^n \times (0, \tau)$, therefore
\[
\log \hat \rho(x, t) \le C_0\left(1 + |\log \eps| +  |\log t|  + \left| \log \frac{t^2}{\eps} \right| \right).
\]
Then, by the conservation of mass
  \begin{equation}\label{LlogLbound}
  \begin{split}
\int_{\T^n}\hat \rho(y,t) |\log \hat \rho(y, t)| \dd y &= - \int_{\T^n\cap \{\hat \rho \leq 1\}}\hat \rho(y,t) \log \hat \rho(y, t) \dd y+\int_{\T^n\cap \{\hat \rho > 1\}}\hat \rho(y,t) \log \hat \rho(y, t) \dd y \\
&\le 1 + C_0\left(1 + |\log \eps| +  |\log t|  + \left| \log \frac{t^2}{\eps} \right| \right) \int_{\T^n\cap \{\hat \rho > 1\}}\hat \rho(y,t)  \dd y  \\
& \le C_1\left(1 + |\log \eps| +  |\log t|  + \left| \log \frac{t^2}{\eps} \right| \right).
\end{split}
\end{equation}


Let us now test the equation for $\hat \rho$ by $\log \hat \rho$ and integrate to obtain, for $t \in (0,\tau)$,
\[
\frac{\dd}{\dd t}\int_{\T^n} \hat \rho(t) \log \hat \rho(t)\dd x + \eps \int_{\T^n} \frac{|D \hat \rho(t)|^2}{\hat \rho(t)}\dd x = \int_{\T^n} (\hat b D\hat \rho)(t)\dd x,
\]
hence by Young's inequality and the conservation of mass
\[
\frac{\dd}{\dd t}\int_{\T^n} \hat \rho(t) \log \hat \rho(t)\dd x + \frac \eps 2 \int_{\T^n} \frac{|D \hat \rho(t)|^2}{\hat \rho(t)}\dd x \le \frac{1}{2\eps}\int_{\T^n} |\hat b|^2 \hat \rho(t)\dd x
\le \frac{K^2}{2 \eps}.
\]
Integrating on intervals $(t_1,t_2) \subset (0, \tau)$ and plugging in \eqref{LlogLbound} gives then
\begin{equation}\label{Fisherbound}
\eps \int_{t_1}^{t_2} \int_{\T^n} \frac{|D \hat \rho(t)|^2}{\hat \rho(t)} \dd t\le C_2\left(1 + |\log \eps| +  |\log t_1| +  |\log t_2| + \left| \log \frac{t_1^2}{\eps} \right| + \left| \log \frac{t_2^2}{\eps} \right| + \frac{t_2-t_1}{\eps} \right).
\end{equation}

Choose now $t_1 = \eps$ and $t_2 = 2\eps$. By the Mean Value Theorem there exists $\hat t \in [\eps, 2\eps]$ such that
\[
\int_{\T^n} \frac{|D \hat \rho(\hat t)|^2}{\hat \rho(\hat t)}\dd x= \frac1\eps \int_{\eps}^{2\eps} \int_{\T^n} \frac{|D \hat \rho(t)|^2}{\hat \rho(t)}\dd x \dd t \le \frac{C_3}{\eps^2}\left(1 + |\log \eps|\right).
\]
Thus, by H\"older's inequality we get
\[
\int_{\T^n} |D \hat \rho(\hat t)|\dd x \le \left(\int_{\T^n} \frac{|D \hat \rho(\hat t)|^2}{\hat \rho(\hat t)}\dd x\right)^{1/2}\left(\int_{\T^n} {\hat \rho(\hat t)}\dd x\right)^{1/2} \le \frac{C_3^{1/2}}{\eps}\left(1 + |\log \eps|\right)^{1/2}.
\]
Similarly, by choosing $t_1 = T- \tau - 2\eps$ and $t_2 = T- \tau - \eps$ we have the existence of $\tilde t \in [T- \tau - 2\eps, T- \tau - \eps]$ such that
\[
\int_{\T^n} |D \hat \rho(\tilde t)|\dd x \le \frac{C_4^{1/2}}{\eps}\left(1 + |\log \eps|\right)^{1/2}.
\]
Going back to the original time variable we obtain the assertion.
\end{proof}

In what follows, we will use again the adjoint problem for $\rho= \rho_{\eps, \bar x, \tau}$, where $\bar x \in \T^n$, $\tau \in [0,T)$ 
\[
\begin{cases} 
\partial_t \rho - \varepsilon \Delta \rho - { \rm div}(D_pH(Du_\eps) \rho) = 0 & \text{in } Q_{\tau,T}, \\
\rho(\tau) = \delta_{\bar x}, & \bar x \in \T^n.
\end{cases}
\]

\begin{lemma}\label{crosslip} Assume that $0 < \eps \le 1$ and $\tau \le T-4\eps$. Then, there exists $C$ depending on $n, \|Du_\eps(T)\|_{L^\infty_x}, T$, and $t_1,t_2$ satisfying
\[
t_1 \in [\tau+\eps, \tau + 2\eps], \qquad t_2 \in [T-2\eps, T-\eps]
\]
such that
\[
 \int_{t_1}^{t_2}\int_{\T^n} (t-\tau)(T-t) |D^2u_\eps|^2\rho\,\dd x \dd t \le C(1 + |\log \eps|).
\]
\end{lemma}

\begin{proof}
We find, by uniform convexity of $H$, the following inequality solved by the function $z(x,t)= \chi(t) \Delta u_\eps (x,t)$, where $\chi(t) = (t-\tau)(T-t)$ : 
\[
-\partial_t z-\eps\Delta z+\theta \chi |D^2u_\eps|^2+DH(Du_\eps)\cdot Dz\leq -\chi' \Delta u_\eps  + \chi \Delta f\text{ in }Q_{\tau,T}.
\]

Let now $t_1,t_2$ be as in Lemma \ref{lemDrho}, applied with $b = -D_pH(Du_\eps)$, that gives
 \begin{equation}\label{drhoest2}
 \int_{\T^n} |D \rho(t_i)| \dd x \le \frac{C}{\eps}\left(1 + |\log \eps|\right)^{1/2}, \qquad i=1,2.
 \end{equation}
 for some
 \begin{equation}\label{tiwhere}
t_1 \in [\tau+\eps, \tau + 2\eps], \qquad t_2 \in [T-2\eps, T-\eps].
\end{equation}

By duality between $\rho$ and $z$ and integrating in $\T^n\times(t_1, t_2) $ we have
\begin{multline}\label{eqzz1}
\theta\int_{t_1}^{t_2}\int_{\T^n}\chi |D^2u_\eps|^2\rho\,\dd x \dd t\leq \int_{\T^n}z(t_2)\rho(t_2) \dd x- \int_{\T^n}z(t_1)\rho(t_1) \dd x - \int_{t_1}^{t_2}\int_{\T^n} \chi' \Delta u_\eps \rho \dd x \dd t + \\ + \int_{t_1}^{t_2}\int_{\T^n} \chi \Delta f \rho \dd x \dd t.
\end{multline}
Note that $(\Delta u_\eps)^2\leq n|D^2u_\eps|^2$, hence by Young's inequality we get
\begin{multline*}
- \int_{t_1}^{t_2}\int_{\T^n} \chi' \Delta u_\eps \rho \dd x \dd t \le
\frac \theta 2 \int_{t_1}^{t_2}\int_{\T^n} \chi |D^2 u_\eps|^2 \rho \dd x \dd t + \frac {n}{2\theta} \int_{t_1}^{t_2} \frac{(\chi' )^2}{\chi} \left(\int_{\T^n} \rho \dd x \right) \dd t \le 
\\ \frac \theta 2 \int_{t_1}^{t_2}\int_{\T^n} \chi |D^2 u_\eps|^2 \rho \dd x \dd t + \\ + \frac {C}{2\theta} \left(4(t_2-t_1) + (T-\tau)(|\log(t_2-\tau)| + |\log(t_1-\tau)| + |\log(T-t_2)| + |\log(T-t_1)| ) \right).
\end{multline*}
On the other hand, for $i=1,2$, integrating by parts yields, together with Lemma \ref{lemDrho} and the Lipschitz estimates of Lemma \ref{lipbound}
\[
\left| \int_{\T^n}z(t_i)\rho(t_i) \dd x \right| \le \chi(t_i)  \int_{\T^n}|Du_\eps(t_i)| \, |D \rho(t_i)| \dd x \le (t_i - \tau)(T-t_i) \frac{C\left(1 + |\log \eps|\right)^{1/2}}{\eps}.
\]

Plugging the previous inequality into \eqref{eqzz1}, and using also \eqref{tiwhere}, we conclude that
\[
\frac\theta 2 \int_{t_1}^{t_2}\int_{\T^n}\chi |D^2u_\eps|^2\rho\,\dd x \dd t \le C_5\left((1 + |\log \eps|)^{1/2} + 1 + |\log \eps|\right)
\]
for some $C_5$ depending on $n, \|Du_T\|_{L^\infty_x}, T,  \|(\Delta f)^+\|_{L^\infty_{x,t}}$.
\end{proof}

\begin{proof}[Proof of Theorem \ref{main2}]
As in Theorem \ref{main1} we have to estimate $\int_\tau^T\int_{\T^n}\Delta u_\eps\rho\,\dd x\dd t$. Let $t_1$ and $t_2$ be as in Lemma \ref{lemDrho} and $\tau \le T-4\eps$. Then,
\begin{multline*}
\left|\int_\tau^T\int_{\T^n}\Delta u_\eps\rho\,\dd x\dd t \right| \leq \sqrt{n}\int_\tau^T\int_{\T^n}|D^2 u_\eps|\rho\,\dd x\dd t \\
= \sqrt{n}\int_\tau^{t_1}\int_{\T^n}|D^2 u_\eps|\rho\,\dd x\dd t + \sqrt{n}\int_{t_2}^{T}\int_{\T^n}|D^2 u_\eps|\rho\,\dd x\dd t + \sqrt{n}\int_{t_1}^{t_2}\int_{\T^n}|D^2 u_\eps|\rho\,\dd x\dd t.
\end{multline*}

First, since $t_1-\tau \le 2\eps$ and $T-t_2 \le 2\eps$, by Lemma \ref{lipbound}
\[
\int_\tau^{t_1}\int_{\T^n}|D^2 u_\eps|\rho\,\dd x\dd t 
\le \left(\int_\tau^{t_1}\int_{\T^n}|D^2 u_\eps|^2\rho\,\dd x\dd t\right)^{1/2}(t_1-\tau)^{1/2} \\ \leq \left( \frac{C}{\eps}\right)^{1/2}(2\eps)^{1/2}=(2C)^\frac12,
\]
and similarly
\[
\int_{t_2}^{T}\int_{\T^n}|D^2 u_\eps|\rho\,\dd x\dd t \le (2C)^{1/2}.
\]

On the other hand, by Lemma \ref{crosslip},
\begin{align*}
\left(\int_{t_1}^{t_2}\int_{\T^n}|D^2 u_\eps|\rho\,\dd x\dd t \right)^2 &= \left(\int_{t_1}^{t_2}\int_{\T^n}\frac{(t-\tau)^{1/2}(T-t)^{1/2}}{(t-\tau)^{1/2}(T-t)^{1/2}}|D^2 u_\eps|\rho\,\dd x\dd t \right)^2\\ 
&\le\left(\int_{t_1}^{t_2}\int_{\T^n} (t-\tau)(T-t) |D^2 u_\eps|^2\rho\,\dd x\dd t \right)\left(\int_{t_1}^{t_2} \frac1{(t-\tau)(T-t)} \dd t\right) \\ 
&\le C_6(1 + |\log \eps|)\frac{|\log|t_1-T|| +|\log|t_1-\tau||+|\log|t_2-T|| + |\log|t_2-\tau||}{T-\tau} \\ 
&\le C_7(1 + |\log \eps|)\frac{|\log \eps|}{T-\tau}.
\end{align*}
Therefore,
\[
\left|\int_\tau^T\int_{\T^n}\Delta u_\eps\rho\dd x\dd t \right| \leq C_8\frac{1+|\log \eps|}{(T-\tau)^{1/2}}.
\]
The rate of convergence 
\[
\|(u_\eps-u)(\tau)\|_{L^\infty(\T^n)}\leq \frac{C_9}{(T-\tau)^{1/2}}\eps|\log\eps|
\]
is straightforward from Remark \ref{concl}.
\end{proof}

\section{Rate of convergence for Hamiltonians with superquadratic growth and semiconcave data}\label{sec;strictly}
In this section we consider strictly convex Hamiltonians with superquadratic growth. 
We will need to prove a refined version of \eqref{secondord} under \eqref{H2}: this will be performed through a Bernstein-type argument.

\begin{thm}\label{main3}
Let $H$ be such that \eqref{H2} holds. Then, if $u_\eps$ solves \eqref{HJintroVisc} and $u$ solves the first-order equation \eqref{HJintro}, we have for some $\beta_\gamma\in\left(\frac12,1\right)$
\[
-C\eps^{\beta}\leq u_\eps-u\leq (T\|(\Delta u_T)^+\|_{L^\infty_x}+\|(\Delta f)^+\|_{L^1_t(L^\infty_x)})\eps \quad \text{ for all } \beta<\beta_\gamma.
\]
where $C$ depends on $n,\theta,\gamma,\beta,\|(\Delta u_T)^+\|_{L^\infty_x}, T, \|(\Delta f)^+\|_{L^\infty_{x,t}}, \|Du_T\|_{L^\infty},\|Df\|_{L^\infty_{x,t}}$ and $\beta_\gamma$ depends on $\gamma$.
\end{thm}

The value of $\beta_\gamma$ is explicit, see \eqref{gammaval} below.

\begin{proof}
The proof of the bound from above is the same of Theorem \ref{main1}, owing to the semiconcavity estimate in \cite[Proposition 3.7 and Remark 3.8]{CGsima}. We prove the bound from below following similar steps. \\

\textit{Step 1}. We prove that for $m$ small enough and $\eta < 1/2$ it holds
\begin{equation}\label{so2}
\iint_{Q_\tau\cap \{|Du|^2\leq m\}}|D^2u_\eps|^2\rho\dd x \dd t\leq \frac{C_\eta m^\eta}{\eps}.
\end{equation}
Let us first observe that for $\tilde w=\varphi(|Du|^2)$, $\varphi\in C^2(\R)$, we have
\[
\Delta \tilde w=4\sum_i \varphi''(Du_{x_i}\cdot Du)^2+2\varphi'D\Delta u\cdot Du+2\varphi'|D^2u|^2,
\]
which leads to the evolution PDE solved by $\tilde w$
\begin{equation}\label{tildew}
-\partial_t \tilde w-\eps\Delta \tilde w+2\eps\varphi'|D^2u|^2+D_pH(Du)\cdot D\tilde w=-4\eps \sum_i \varphi''(Du_{x_i}\cdot Du)^2+2\varphi'Df\cdot Du.
\end{equation}
Notice that if $\varphi$ is concave, we get
\[
-\partial_t \tilde w-\eps\Delta \tilde w+2\eps\varphi'(|Du|^2)|D^2u|^2+D_pH(Du)\cdot D\tilde w \le -4\eps \varphi''( |Du|^2) |Du|^2 |D^2u|^2+2\varphi'Df\cdot Du.
\]
Let now $ 0 < \delta < 1/2$, and choose 
\[
\varphi(z) = \varphi_m(z) = 
\begin{cases}
z & z\in[0,m] \\
2(\delta+1) \sqrt{mz} - \delta z -\delta m-m & z \in [m,\kappa m ], \qquad \kappa = \frac{(1+\delta)^2}{\delta^2} >1\\
m\frac{1+\delta}{\delta} & z \in [\kappa m, +\infty).
\end{cases}
\]
One may verify with a direct computation that $\varphi$ is nonnegative, increasing, concave and $C^1$. Moreover,
\begin{align*}
&\varphi' = 1, \qquad \varphi'' = 0 \quad \text{on $[0,m)$} \\
&\varphi'(z) + 2 \varphi''(z) z = -\delta \quad \text{on $(m,\kappa m)$} \\
&\varphi' = \varphi'' = 0 \quad \text{elsewhere.}
\end{align*}
Using these properties, and testing the PDE \eqref{tildew} by $\rho$ and integrating  
we conclude\footnote{Since $\varphi_m$ is not $C^2$, one should argue along a sequence of smooth, concave $\varphi_{m,n}$ that coincide with $\varphi$ on $[0,m] \cup [\kappa m, +\infty)$. Then, the inequality below is obtained in the limit $\varphi_{m,n} \to \varphi_{m}$.}
\begin{multline*}
2\eps\iint_{Q_\tau}\varphi' |D^2u|^2\rho\dd x \dd t \le \int_{\T^n} \tilde w \rho(x,t)\dd x\Big|^{t=T}_{t=\tau}- 4\eps\iint_{ Q_\tau }\varphi'' |Du|^2 |D^2u|^2\rho\dd x \dd t\\
+2\|Df\|_{L^\infty(Q_\tau)}\iint_{Q_\tau}\varphi'|Du|\dd x\dd t,
\end{multline*} 
hence, using $\varphi'(z) + 2 \varphi''(z) z = -\delta$ and $-\delta\leq \varphi'\leq 1$ on $(m,\kappa m)$, we have
\begin{multline*}
\eps\iint_{Q_\tau\cap \{|Du|^2\leq m\}}|D^2u|^2\rho\dd x \dd t\leq m\frac{1+\delta}{\delta}+ 2\delta \eps\iint_{ Q_\tau\cap \{m < |Du|^2 < \kappa m\} } |D^2u|^2\rho\dd x \dd t\\
+2\|Df\|_{L^\infty(Q_\tau)}\iint_{Q_\tau\cap \{|Du|^2 < \kappa m\}}\varphi'|Du|\dd x\dd t.
\end{multline*}
We now proceed by a hole-filling type technique, that is, we ``fill the hole'' by adding to both sides $2\delta \eps\iint_{ Q_\tau\cap \{|Du|^2 \le  m\} } |D^2u|^2\rho$ and obtain
\[
\eps\iint_{Q_\tau\cap \{|Du|^2\leq m\}}|D^2u|^2\rho\dd x \dd t\leq m c_\delta + \frac{2\delta}{1+2\delta} \eps\iint_{ Q_\tau\cap \{ |Du|^2 < \kappa m\} } |D^2u|^2\rho\dd x \dd t+\frac{c_f}{1+2\delta}\sqrt{\kappa m},
\]
where $c_\delta > 0$ depends on $\delta$ only. If we now let, for $m \in (0,1)$,
\[
h(m) := \frac{\eps\iint_{Q_\tau\cap \{|Du|^2\leq m\}} |D^2u|^2\rho\dd x \dd t}{m^\eta} ,
\]
the previous inequality reads, after dividing by $m^\eta$ and using that $m<1$ and $\eta<1/2$, as
\begin{equation}\label{etaobstruction}
h(m) \le c_\delta m^{1-\eta} + \frac{2\delta \kappa^\eta }{1+2\delta} h(\kappa m)+c_f\sqrt{\kappa}m^{\frac12-\eta}\leq c_\delta  + \frac{2\delta \kappa^\eta }{1+2\delta} h(\kappa m)+c_f\sqrt{\kappa}.
\end{equation}
Recalling the definition of $\kappa$ above, on one hand we have $\sqrt{\kappa}=\frac{1+\delta}{\delta}$, and, on the other hand, we can pick $\delta$ small enough so that
\[
\frac{2\delta  }{1+2\delta}\left(\frac{1+\delta}{\delta}\right)^{2\eta} \le \frac 12.
\]
Consider now the sequence $m_j = 1/\kappa^j$, $j = 1, 2, \ldots$ and evaluate the previous inequality on $m_j$ to get
\[
h(m_j) \le \tilde c_\delta + \frac12 h(m_{j-1}),
\]
which yields by induction that $h(m_j) \le 2 \tilde c_\delta + \frac12h(m_1)$ for all $j \ge 2$. Since $h$ is increasing and $h(m_1)$ is bounded by Lemma \ref{lipbound}, we get the claim.

\medskip

\textit{Step 2}. We now proceed by differentiating twice the equation  and use \eqref{H2} to find the following inequality solved by the function $z(x,t)=(t-\tau)^\alpha \Delta u_\eps$, ($\alpha>0$ to be chosen later)
\begin{multline*}
-\partial_t z-\eps\Delta z+\theta (t-\tau)^\alpha |Du_\eps|^{\gamma-2}  |D^2u_\eps|^2+D_pH(Du_\eps)\cdot Dz\\
\leq -\alpha (t-\tau)^{\alpha-1}\Delta u_\eps+ \Delta f\text{ in }Q_T.
\end{multline*}

By duality and integrating in $\T^n\times(\tau,T)$ we have
\begin{multline*}
\underbrace{\int_{\T^n}z(\tau)\rho(\tau)\,dx}_{=0}+\theta\int_\tau^T\int_{\T^n}(t-\tau)^\alpha|Du_\eps|^{\gamma-2}|D^2u_\eps|^2\rho\dd x \dd t= \underbrace{\int_{\T^n}z(T)\rho(T)\,dx}_{\leq (T-\tau)^\alpha\|(\Delta u_T)^+\|_{L^\infty(\T^n)}}\\
-\alpha \int_\tau^T\int_{\T^n} (t-\tau)^{\alpha-1}\Delta u_\eps(x,t)\rho\dd x \dd t+\int_\tau^T\int_{\T^n} (t-\tau)^{\alpha}\Delta f(x,t) \rho\dd x \dd t.
\end{multline*}
The last integral can be bounded as in Theorem \ref{main1}. We split the first integral in two regions (where the gradient is small and large respectively), apply the H\"older inequality and choose $\alpha>1$ as follows
\begin{align*}
-\alpha& \int_\tau^T\int_{\T^n} (t-\tau)^{\alpha-1}\Delta u_\eps\rho\dd x \dd t=-\alpha\iint_{\{|Du|^2\le m\}}(t-\tau)^{\alpha-1}\Delta u_\eps\rho\dd x \dd t\\
&-\alpha\iint_{\{|Du|^2>m\}}(t-\tau)^{\alpha-1}\Delta u_\eps\rho\dd x \dd t\\
&\leq \alpha\sqrt{n}\left(\iint_{\{|Du|^2\le m\}}|D^2u_\eps|^2\rho\dd x \dd t\right)^\frac12\left(\iint_{Q_\tau}(t-\tau)^{2\alpha-2}\rho\dd x \dd t\right)^\frac12\\
&\qquad +\alpha\sqrt{n}\left(\iint_{\{|Du|^2>m\}}(t-\tau)^{\alpha}|D^2u_\eps|^2\rho\dd x \dd t\right)^\frac12\left(\iint_{Q_\tau}(t-\tau)^{\alpha-2}\rho\dd x \dd t\right)^\frac12\\
&\leq C_{\alpha,T}\sqrt{C_\eta}\sqrt{\frac{m^\eta}{\eps}}+\frac{\theta}{2}\iint_{Q_\tau}(t-\tau)^\alpha|Du_\eps|^{\gamma-2}|D^2u_\eps|^2\rho\dd x \dd t+\frac{C_{\alpha,T}}{m^{\frac{\gamma-2}{2}}}.
\end{align*}
This implies
\[
\int_\tau^T\int_{\T^n}(t-\tau)^\alpha|Du_\eps|^{\gamma-2}|D^2u_\eps|^2\rho\dd x \dd t\leq \tilde K_1\left(\sqrt{\frac{m^\eta}{\eps}}+\frac{1}{m^{\frac{\gamma-2}{2}}}+\|(\Delta f)^+\|_{L^\infty(Q_T)}\right).
\]
Choosing $m = \eps^{\frac{1}{\eta + \gamma -2}}$ we get the following integral bound on superlevel sets $\{|Du|^2>m\}$
\[
\iint_{Q_\tau\cap \{|Du|^2>m\}}(t-\tau)^\alpha|D^2u_\eps|^2\rho\dd x \dd t\leq  \tilde K_2\left(\eps^{\frac{2-\gamma}{2(\gamma+\eta-2)}}+1\right).
\]
We can now conclude the proof. By the foregoing estimates and Lemma \ref{lipbound}
\begin{align*}
\int_\tau^T\int_{\T^n}\Delta u_\eps\rho\dd x \dd t&\leq \sqrt{n}\int_{t_0}^{T}\int_{\T^n}(t-\tau)^{\alpha/2}|D^2 u_\eps|(t-\tau)^{-\alpha/2}\rho\dd x \dd t+\sqrt{n}\int_{\tau}^{t_0}\int_{\T^n}|D^2 u_\eps|\rho\dd x \dd t\\
&\leq \sqrt{n}\left(\int_{t_0}^{T}\int_{\T^n}(t-\tau)^{-\alpha}\rho\dd x \dd t\right)^\frac12\left(\int_{t_0}^{T}\int_{\T^n}|D^2 u_\eps|^2(t-\tau)^{\alpha}\rho\dd x \dd t\right)^\frac12\\
&\quad +\sqrt{n}\frac{C_L}{\sqrt{\eps}}\sqrt{t_0-\tau}
\leq \frac{\sqrt{n}}{|\alpha-1|}((t_0-\tau)^{1-\alpha}+(T-\tau)^{1-\alpha})^\frac12\times\\
&\left(\iint_{\{|Du|^2 \le m\}}|D^2 u_\eps|^2(t-\tau)^{\alpha}\rho\dd x \dd t+\iint_{Q_\tau\cap \{|Du|^2> m\}}|D^2 u_\eps|^2(t-\tau)^{\alpha}\rho\dd x \dd t\right)^\frac12\\
& \quad+\sqrt{n}\frac{C_L}{\sqrt{\eps}}\sqrt{t_0-\tau}\\
&\leq  \tilde K_3\sqrt{n}((t_0-\tau)^{1-\alpha}+(T-\tau)^{1-\alpha})^\frac12\left(\frac{m^\eta}{\eps}+ \eps^{\frac{2-\gamma}{2(\gamma+\eta-2)}} +1\right)^\frac12\\
&\quad +\sqrt{n}\frac{C_L}{\sqrt{\eps}}\sqrt{t_0-\tau}.
\end{align*}
With the choice of $m = \eps^{\frac{1}{\eta + \gamma -2}}$ as before we get
\[
\int_\tau^T\int_{\T^n}\Delta u_\eps\rho\dd x \dd t\leq \tilde K_4 \left((t_0-\tau)^{1-\alpha}\eps^{\frac{2-\gamma}{2(\gamma+\eta-2)}}+1 \right)+\sqrt{n}\frac{C_L}{\sqrt{\eps}}\sqrt{t_0-\tau}.
\]
If $\tau + \eps  < T$ we choose $t_0=\tau+\eps$ and $\alpha$ arbitrarily close to 1, to get the following result integrating with respect to $\eps$ via Remark \ref{concl}
\[
\|u_\eps-u\|_{L^\infty(Q_T)}\leq  \tilde K_5\eps^{\beta},\ \beta<\frac{\gamma+2\eta-2}{2(\gamma+\eta-2)}.
\]
Then, the claim follows by choosing 
\begin{equation}\label{gammaval}
\beta_\gamma = \sup_{\eta \in (0,1/2)} \frac{\gamma+2\eta-2}{2(\gamma+\eta-2)} = \frac{\gamma-1}{2\gamma-3}.
\end{equation}

\end{proof}

\begin{rem}[On the optimal convergence rate for strictly convex Hamiltonians] \label{superquadrem} Note first that
\begin{align*}
& \beta_\gamma = \frac{\gamma-1}{2\gamma-3} \to 1 \qquad \text{as $\gamma \to 2$},\\
& \beta_\gamma = \frac{\gamma-1}{2\gamma-3} \to \frac 12 \qquad \text{as $\gamma \to \infty$}.
\end{align*}
These limit rates are coherent with the rates obtained in the quadratic case and in the locally Lipschitz one. Our results are in the direction suggested by the numerical experiments in \cite{Sprekeleretal}, that show a rate of order $\eps^\beta$, $\beta\in(1/2,1)$ in the case of strictly convex Hamiltonians. Example 5 therein for instance indicates the order $\mathcal{O}(\eps^\frac23)$, when $H(p)=\frac14|p|^4$ and $f = 0$. Our endpoint rate $\beta_\gamma$, once specialized to $\gamma = 4$ becomes $3/5$, that is slightly worse than $2/3$. We do not know what could be the optimal convergence rate for $H(p)=\frac14|p|^4$, and we cannot exclude that it may depend on the regularity of $f$. 

Note that, in our argument, the rate could be improved if we could allow $\eta$ to vary in the wider range $(0,1)$; for $\eta \to 1$, the exponent $2/3$ could be reached for quartic Hamiltonians: for general $\gamma > 2$, it would be $\frac{\gamma}{2(\gamma-1)}$. Is this the optimal one? Here, we can allow for $\eta < 1/2$ (see for example equation \eqref{etaobstruction}), but that step of the proof ``sees'' only the Lipschitz property of $f$. We do not know now how to improve this step using the information that $f$ is more regular.

\end{rem}

\begin{rem}[On the convergence rate for some Mean Field Control problems]\label{mfc}

The issue of the convergence problem in Mean Field Control has been in the last few years an active area of research, see for instance \cite{DaudinDelarueJackson} and references therein. In general, it amounts to study the convergence of (symmetric) value functions of some control problems as the dimension $n$ of the Euclidean space increases, towards a limit value function, which satisfies an equation in the space of probability measures. In some special cases, the problem can be recast into a vanishing viscosity approximation of (finite-dimensional) Hamilton-Jacobi equations, where in particular the viscosity satisfies $\eps = \eps_n = 1/n$ (see in particular \cite[Proposition 2.10]{DaudinDelarueJackson}). The convergence rate in Mean Field Control is known in general (under Lipschitz regularity assumptions) to be of order $\mathcal{O} (1/\sqrt {n})$; under convexity properties required by Theorems \ref{main1}, \ref{main2} or \ref{main3}, an improved rate of order $1/n^\beta$, $\beta\in(1/2,1)$, or $(1/n)|\log(1/n)|$ shows up in the specific examples described in  \cite{DaudinDelarueJackson} .

\end{rem}


\begin{thebibliography}{QSTY24}

\bibitem[BCC22]{BCC}
P. Bonicatto, G. Ciampa, and G. Crippa.
\newblock On the advection-diffusion equation with rough coefficients: weak
  solutions and vanishing viscosity.
\newblock {\em J. Math. Pures Appl. (9)}, 167:204--224, 2022.

\bibitem[BCD97]{BCD}
M. Bardi and I. Capuzzo-Dolcetta.
\newblock {\em Optimal control and viscosity solutions of
  {H}amilton-{J}acobi-{B}ellman equations}.
\newblock Systems \& Control: Foundations \& Applications. Birkh\"{a}user
  Boston, Inc., Boston, MA, 1997.

\bibitem[BKRS15]{BKRS}
V.~I. Bogachev, N.~V. Krylov, M. R\"{o}ckner, and S.~V.
  Shaposhnikov.
\newblock {\em Fokker-{P}lanck-{K}olmogorov equations}, volume 207 of {\em
  Mathematical Surveys and Monographs}.
\newblock American Mathematical Society, Providence, RI, 2015.


\bibitem[Bre20]{Brenier}
Y.~Brenier.
\newblock Examples of hidden convexity in nonlinear {PDEs}, {L}ecture {N}otes
  at {ETH} {Z}urich.
\newblock hal-02928398, 2020.

\bibitem[Cal18]{Calder}
J.~Calder.
\newblock {L}ecture notes on viscosity solutions.
\newblock available at
  \url{https://www-users.cse.umn.edu/~jwcalder/viscosity_solutions.pdf}, 2018.


\bibitem[CCM11]{CCM}
F. Camilli, A. Cesaroni, and C. Marchi.
\newblock Homogenization and vanishing viscosity in fully nonlinear elliptic
  equations: rate of convergence estimates.
\newblock {\em Adv. Nonlinear Stud.}, 11(2):405--428, 2011.

\bibitem[CD25]{ChaintronDaudin}
L.-P. Chaintron and S.~Daudin.
\newblock Optimal rate of convergence in the vanishing viscosity for uniformly convex
  {H}amilton-{J}acobi equations, 
\newblock arXiv: 2502.09103, 2025.

\bibitem[CG19]{CGsima}
M. Cirant and A. Goffi.
\newblock On the existence and uniqueness of solutions to time-dependent
  fractional {MFG}.
\newblock {\em SIAM J. Math. Anal.}, 51(2):913--954, 2019.

\bibitem[CM24]{alpar}
M. Cirant and A. R. M\'esz\'aros.
\newblock Long Time Behavior and Stabilization for Displacement Monotone Mean Field Games.
\newblock arXiv:2412.14903, 2024.

\bibitem[CGM23]{CGM}
F. Camilli, A. Goffi, and C. Mendico.
\newblock Quantitative and qualitative properties for {Hamilton}-{J}acobi
  {PDE}s via the nonlinear adjoint method, to appear in Ann. Sc. Norm. Super.
  Pisa, Cl. Sci., 2023.
  
 \bibitem[CGMT15]{CGMT}
F.~Cagnetti, D.~Gomes, H.~Mitake, and H.~V. Tran.
\newblock A new method for large time behavior of degenerate viscous
  {H}amilton-{J}acobi equations with convex {H}amiltonians.
\newblock {\em Ann. Inst. H. Poincar\'{e} Anal. Non Lin\'{e}aire},
  32(1):183--200, 2015. 
 
  \bibitem[CGT13]{CGT}
F. Cagnetti, D. Gomes, and H. V.Tran.
\newblock Convergence of a semi-discretization scheme for the {Hamilton}-{Jacobi} equation: a new approach with the adjoint method, Appl. Numer. Math. 73:2--15, 2013.

\bibitem[Con24]{conforti}
G.~Conforti.
\newblock Weak semiconvexity estimates for Schrödinger potentials and logarithmic Sobolev inequality for Schrödinger bridges.
  \newblock {\em Probability Theory and Related Fields }, 189, 1045–-1071, 2024.

\bibitem[CL84]{CL84}
M.~G. Crandall and P.-L. Lions.
\newblock Two approximations of solutions of {H}amilton-{J}acobi equations.
\newblock {\em Math. Comp.}, 43(167):1--19, 1984.

\bibitem[DDJ24]{DaudinDelarueJackson}
S. Daudin, F. Delarue, and J. Jackson.
\newblock On the optimal rate for the convergence problem in mean field
  control.
\newblock {\em J. Funct. Anal.}, 287(12):94, 2024.
\newblock Id/No 110660.

\bibitem[DI06]{DroniouImbert}
J. Droniou and C. Imbert.
\newblock Fractal first-order partial differential equations.
\newblock {\em Arch. Ration. Mech. Anal.}, 182(2):299--331, 2006.

\bibitem[Eva10]{EvansARMA}
L.~C. Evans.
\newblock Adjoint and compensated compactness methods for {H}amilton-{J}acobi
  {PDE}.
\newblock {\em Arch. Ration. Mech. Anal.}, 197(3):1053--1088, 2010.

\bibitem[Fle64a]{FlemingPrinceton}
W.~H. Fleming.
\newblock The convergence problem for differential games, {II}.
\newblock {\em Ann. Math. Stud.}, 52:195--210, 1964.

\bibitem[Fle64b]{FlemingJMM}
W.~H. Fleming.
\newblock The Cauchy problem for degenerate parabolic equations.
\newblock {\em J. Math. and Mech.}, 13:987--1008, 1964.

\bibitem[FS86]{FlemingSouganidis}
W.~H. Fleming and P.~E. Souganidis.
\newblock Asymptotic series and the method of vanishing viscosity.
\newblock {\em Indiana Univ. Math. J.}, 35(2):425--447, 1986.

\bibitem[Gof24]{Goffi}
A. Goffi.
\newblock Remarks on the rate of convergence of the vanishing viscosity process
  of {H}amilton-{J}acobi equations.
\newblock arXiv:2412.15651, 2024.




\bibitem[Kru66]{K66d}
S.~N. Kruzhkov.
\newblock Solutions of first-order nonlinear equations.
\newblock {\em Sov. Math., Dokl.}, 7:376--379, 1966.

\bibitem[Kru67]{K67II}
S.~N. Kruzhkov.
\newblock Generalized solutions of nonlinear equations of the first order with
  several independent variables. {II}.
\newblock {\em Mat. Sb. (N.S.)}, 72 (114):108--134, 1967.

\bibitem[IK91]{IshiiKoike}
H. Ishii and S. Koike.
\newblock {\em Remarks on elliptic singular perturbation problems}.
\newblock {\em Appl. Math. Optim.} 23:1--15, 1991.


\bibitem[LBL19]{LBL}
C. Le~Bris and P.-L. Lions.
\newblock {\em Parabolic equations with irregular data and related
  issues---applications to stochastic differential equations}, volume~4 of {\em
  De Gruyter Series in Applied and Numerical Mathematics}.
\newblock De Gruyter, Berlin, 2019.

\bibitem[Lio82]{L82Book}
P.-L. Lions.
\newblock {\em Generalized solutions of {H}amilton-{J}acobi equations},
  volume~69 of {\em Research Notes in Mathematics}.
\newblock Pitman (Advanced Publishing Program), Boston, Mass.-London, 1982.

\bibitem[Lio85]{Lions85aa}
P.-L. Lions.
\newblock Regularizing effects for first-order {H}amilton-{J}acobi equations.
\newblock {\em Applicable Anal.}, 20(3-4):283--307, 1985.

\bibitem[LMT17]{LMT}
N.Q. Le, H.~Mitake and H.~V. Tran.
\newblock {\em Dynamical and geometric aspects of {H}amilton-{J}acobi and
  linearized {M}onge-{A}mp\`ere equations---{VIASM} 2016}, volume 2183 of {\em
  Lecture Notes in Math.}, pages 125--228. Springer, Cham, 2017.

\bibitem[LSU68]{LSU}
O.~A. Lady\v{z}enskaja, V.~A. Solonnikov, and N.~N. Ural'~ceva.
\newblock {\em Linear and quasilinear equations of parabolic type}.
\newblock Translations of Mathematical Monographs, Vol. 23. American
  Mathematical Society, Providence, RI, 1968.
\newblock Translated from the Russian by S. Smith.

\bibitem[LT01]{LinTadmor}
C.-T. Lin and E. Tadmor.
\newblock {$L^1$}-stability and error estimates for approximate
  {H}amilton-{J}acobi solutions.
\newblock {\em Numer. Math.}, 87(4):701--735, 2001.

\bibitem[NFSS22]{NSS}
V. Navarro-Fern\'{a}ndez, A. Schlichting, and C. Seis.
\newblock Optimal stability estimates and a new uniqueness result for
  advection-diffusion equations.
\newblock {\em Pure Appl. Anal.}, 4(3):571--596, 2022.

\bibitem[NT92]{NT}
H. Nessyahu and E. Tadmor.
\newblock The convergence rate of approximate solutions for nonlinear scalar
  conservation laws.
\newblock {\em SIAM J. Numer. Anal.}, 29(6):1505--1519, 1992.

\bibitem[QSTY24]{Sprekeleretal}
J. Qian, T. Sprekeler, H.~V. Tran, and Y. Yu.
\newblock Optimal rate of convergence in periodic homogenization of viscous
  {Hamilton}-{Jacobi} equations.
\newblock {arXiv}:2402.03091, to appear in
  Multiscale Modeling and Simulation, 2024.

\bibitem[Tra11]{TranCalcVarPDE}
H.~V. Tran.
\newblock Adjoint methods for static {H}amilton-{J}acobi equations.
\newblock {\em Calc. Var. Partial Differential Equations}, 41(3-4):301--319,
  2011.

\bibitem[Tra21]{TranBook}
H.~V. Tran.
\newblock {\em Hamilton-Jacobi equations: {T}heory and {A}pplications}.
\newblock AMS Graduate Studies in Mathematics. American Mathematical Society,
  2021.
  
  \bibitem[TT99]{TT}
E.~Tadmor and T.~Tang.
\newblock Pointwise error estimates for scalar conservation laws with piecewise smooth solutions.
\newblock {\em SIAM J. Numer. Anal.}, 36(6):1739--1758, 1999.
  
  \bibitem[TZ23]{TangZhang}
W.~Tang and Y.~P. Zhang.
\newblock The convergence rate of vanishing viscosity approximations for mean
  field games.
\newblock arXiv:2303.14560, 2023.



\bibitem[Wan99]{Wang}
W.-C. Wang.
\newblock On {{\(L^{1}\)}} convergence rate of viscous and numerical
  approximate solutions of genuinely nonlinear scalar conservation laws.
\newblock {\em SIAM J. Math. Anal.}, 30(1):38--52, 1999.

\bibitem[Wid71]{Widman}
K.-O. Widman.
\newblock H{\"o}lder continuity of solutions of elliptic systems.
\newblock {\em Manuscr. Math.}, 5:299--308, 1971.

\end{thebibliography}
\end{document}